\def\draft{n}
\documentclass{amsart} 
\usepackage{latexsym,amssymb,graphicx,amscd,amsmath,url}

\def\printname#1{
	\if\draft y
		\smash{\makebox[0pt]{\hspace{-0.5in}
			\raisebox{8pt}{\tt\tiny #1}}}
	\fi
}

\def\lbl#1{\label{#1}\printname{#1}}

\newtheorem{thm}{Theorem}
\newtheorem{lem}[thm]{Lemma}

\newtheorem{cor}[thm]{Corollary}

\newcommand{\BS}{{\mathcal{S}}}

\newcommand{\BZ}{{\mathbb{Z}}}
\newcommand{\BQ}{{\mathbb{Q}}}

\newcommand{\BF}{{\mathbb{F}}}

\newcommand{\Si}{{\Sigma}}

\newcommand{\fo}{\mathfrak{o}}
\newcommand{\fe}{\mathfrak{e}}

\newcommand{\I}{{\mathrm I}}

\DeclareMathOperator{\Sp}{Sp}

\DeclareMathOperator{\SL}{SL}

\DeclareMathOperator{\SU}{SU}
\DeclareMathOperator{\SO}{SO}
\DeclareMathOperator{\diag}{diag}

\begin{document}

  \title[Dimension formulas for some modular representations]
 {Dimension formulas for some modular representations of the
  symplectic group in the natural characteristic}

\author{ Patrick M. Gilmer}
\address{Department of Mathematics\\
Louisiana State University\\
Baton Rouge, LA 70803\\
USA}
\email{gilmer@math.lsu.edu}
\thanks{The first author was partially supported by  NSF-DMS-0905736 }
\urladdr{\url{www.math.lsu.edu/~gilmer}}

\author{Gregor Masbaum}
\address{Institut de Math{\'e}matiques de Jussieu (UMR 7586 du CNRS)\\
Case 247, 
4 pl. Jussieu, 
75252 Paris Cedex 5,
FRANCE }
\email{masbaum@math.jussieu.fr}
\urladdr{\url{www.math.jussieu.fr/~masbaum}}


\date{November 11, 2011}

 \begin{abstract} We compare the dimensions of the
irreducible  $\Sp(2g,K)$-modules over a field $K$ of characteristic
$p$ constructed by Gow \cite{Go} with 
the dimensions of
the irreducible   $\Sp(2g,\BF_p)$-modules that appear  in the
 first approximation to
representations of mapping class groups of surfaces in Integral
Topological Quantum Field Theory \cite{GM4}.  For this purpose, we derive a 
trigonometric  formula for the dimensions of Gow's representations.
 This formula is equivalent to a special case of a  formula contained in unpublished work of Foulle
\cite{F1,F2}.  Our direct proof is simpler than the proof of Foulle's more
general result, and is modeled on the proof of the  Verlinde formula in TQFT. \end{abstract}

\maketitle

Let $g\geq 1$ be an integer. The  irreducible representation of the
symplectic group $\Sp(2g,{\mathbb C})$ with fundamental highest weight $\omega_k$ $(1\leq
k\leq g)$ can be constructed inside the $k$-th exterior power of the
standard representation of $\Sp(2g,{\mathbb C})$ as the kernel of a certain contraction
operator (see for example Fulton and Harris \cite[\S17]{FH}). The
dimension of this representation is
\begin{equation}\lbl{funk}
\binom{2g}{k} - \binom{2g}{k-2}~.
\end{equation}
 Assume now that $p$ is an odd
prime and $K$ is an algebraically closed field of characteristic $p$.  Premet and Suprunenko \cite{PS} showed that if one performs the same
construction over the field $K$, 
the resulting representation of the symplectic group $\Sp(2g,K)$ is
again irreducible with fundamental highest weight $\omega_k$ provided
$p>g$. But if $p\leq g$, this representation may no longer be
irreducible, although it is still a Weyl module for  $\omega_k$, and
so 
it has a unique irreducible quotient which will be an
irreducible $\Sp(2g,K)$-module with highest weight $\omega_k$ (see for example Humphreys
\cite[\S3.1]{H}).  
Gow
considered precisely this situation in \cite{Go}, except that he did
not assume that $K$ is algebraically closed. For  every integer $g\geq
p-1$,  and for the $p-1$ values of $k$ in the range $[g-p+2,g]$,  he
considered natural subquotients of $\Lambda^kV$, where $V\simeq K^{2g}$
is the standard representation of the symplectic group  
$\Sp(2g,K)$, and showed that these subquotients
are irreducible representations of   
$\Sp(2g,K)$. Following Gow, we denote these representations by 
\begin{equation*}\qquad \qquad \qquad \quad V(g,k) \qquad \quad (g-p+2 \leq k \leq g)~.
\end{equation*} Gow also proved by an explicit computation that if $K$
is algebraically closed, then $V(g,k)$ indeed has
fundamental highest weight $\omega_k$ \cite[Corollary~2.4]{Go}.

In this 
note,  
we are interested in the dimensions of Gow's
representations $V(g,k)$. 
For $g=p-1$ these dimensions are still given by the classical formula
(\ref{funk}). But for $g>p-1$ this is no longer the case, 
and the dimension of $V(g,k)$  
now depends 
on the characteristic $p$ of
the field $K$. Gow gave a recursion formula for 
this dimension,  
 but stopped short of solving the recursion except for
$p\leq 5$. 
We solved Gow's recursion for general primes $p$ and thereby found an
explicit trigonometric formula for the dimension of $V(g,k)$ which we
will state in Theorem~\ref{th1} below. 
 Our proof of this dimension formula is analogous to the proof of the famous
Verlinde formula in Topological Quantum Field Theory (TQFT)
which computes the dimensions of the TQFT vector spaces (see Formula
(\ref{Verl}) below).   
 In fact, we came to 
our  formula while trying to compare the
dimensions of certain irreducible modular representations of
$\Sp(2g,\BF_p)$ arising from our theory of 
Integral TQFT \cite{G,GM,GM4} with Gow's representations
$V(g,k)$. Experimentally, we 
found that for  $p=5$, the dimensions coincide, whereas for $p\geq 7$
they do not.   
The aim of this note is to explain why this is so.

Gow's recursion formula expresses the dimension of $V(g,k)$ as a
matrix coefficient, as follows. Let $C$ be the $(p-1)\times(p-1)$
matrix whose non-zero entries are given by $C_{ii}=2$ for $1\leq i\leq
p-1$ and $C_{j,j-1}=C_{j-1,j}=1$ for 
 $2\leq j\leq p-1$.

\begin{thm}[Gow \cite{Go}] 
Let $p$ be an odd prime and $g\geq p-1$.
 For the $p-1$ values of $k$ in the range $[g-p+2,g]$, one has 
\begin{equation*}\lbl{Gowf}
\dim V(g, k)=  (C^g)_{g+1-k,1}~.
\end{equation*}
\end{thm}

\noindent{\bf Remark.} This formula holds also for $g<p-1$ if we
interpret the left hand side as given by the classical formula
(\ref{funk}) in that case.
\vskip 8pt

We will prove the following:
\begin{thm}\lbl{th1} 
Let $p$ be an odd prime and $g\geq 1$. For $1\leq n\leq p-1$ one has 
$$(C^g)_{n,1}= \frac{2^{2g+1}}{p} \sum_{j=1}^{(p-1)/2} \left(\sin\frac{2\pi
  jn}{p}\right) \left(\sin\frac{2\pi j}{p}\right) \left(
\left(\cos\frac{\pi j}{p}\right)^{2g} - (-1)^{n} \left(\sin\frac{\pi
    j}{p}\right)^{2g}\right)$$
\end{thm}
To get an explicit formula for the dimension of $V(g,k)$ from this
theorem, it suffices to put $n=g+1-k$.

\vskip 8pt
\noindent{\bf Remark.} After completing the first version of this
note, we have learned\footnote{We thank Alexander Premet for telling
  us about Foulle's formula.} that a trigonometric formula for the dimensions of the $V(g,k)$ had been
obtained previously by S.~Foulle.  In fact, 
Foulle in his 2004 thesis \cite{F1} 
(see also his 2005 preprint \cite{F2} which apparently was never published) 
presents 
a more general 
trigonometric dimension formula which works also for fundamental highest
weights 
outside the range $[g-p+2,g]$
considered by Gow.
Foulle's formula is based on Premet and Suprunenko's work
\cite{PS} and is in general much more complicated than the expression in
Theorem~\ref{th1}. But we have checked that Foulle's formula 
when specialized to fundamental highest weights
$\omega_k$, for $k$ 
in the range $[g-p+2,g]$, 
 is equivalent to the formula given in
Theorem~\ref{th1}.   
See \cite[Example~4.6]{F2}.  
Still we believe that our direct proof based on Gow's
recursion formula is of interest, as the proof is analogous
to the proof of the 
Verlinde formula in TQFT and seems much simpler than
Foulle's proof of his more general result.
\vskip 8pt

As already mentioned, our main motivation for proving
Theorem~\ref{th1} was that we need it   
to identify the dimensions of the four
representations 
$V(g,k)$ 
of $\Sp(2g,\BF_p)$ 
in
the case $p=5$  
with the dimensions of 
 the 
irreducible modular
representations constructed in \cite{GM4}.   
The answer will be stated
 in Corollary~\ref{cor1} 
 below. 
First, let 
us
briefly review the construction of these representations.

 Let $\Sigma$ be an
oriented surface of genus $g$ with at most one boundary component. The
orientation-preserving mapping class group of $\Si$ has a natural
surjective homomorphism to the symplectic group $\Sp(2g,\BZ)$, and hence to
$\Sp(2g,\BF_p)$ for every $p$. Using this homomorphism, the modular representations of \cite{GM4} are constructed from representations of the
mapping class group of $\Sigma$ arising in Integral TQFT, as follows.
Let $p\geq 5$ be 
a  
prime, and $2c$ an even
integer with $0\leq 2c\leq p-3$. Integral TQFT \cite{G,GM} provides a
representation of a certain central extension of the mapping class group of $\Si$ on a free
module of finite rank over the ring of cyclotomic integers
$\BZ[\zeta_p]$ (here $\zeta_p$ is a primitive $p$-th root of
unity). These ``quantum" representations are denoted by $\BS(\Si_g(2c))$ in
\cite{GM4}. We think of these representations as an integral version of the Reshetikhin-Turaev TQFT
associated to the Lie group $\SO(3)$, a version of which is obtained if one extends coefficients from $\BZ[\zeta_p]$ to the quotient field
$\BQ(\zeta_p)$. The advantage of having coefficients in $\BZ[\zeta_p]$
is that one can reduce coefficients by dividing out by some ideal in
$\BZ[\zeta_p]$. Recall that $\BZ[\zeta_p]/(1-\zeta_p)$ is the finite field
$\BF_p$. Thus we get a modular representation 
$$F(\Si_g(2c))=\BS(\Si_g(2c))/(1-\zeta_p)\BS(\Si_g(2c))$$ over
$\BF_p$. In some
sense, $F(\Si_g(2c))$ is a first approximation to the quantum
representation $\BS(\Si_g(2c))$. 
While the representations $\BS(\Si_g(2c))$ are irreducible when
considered over the
complex numbers 
(or 
over  quotient fields of $\BZ[\zeta_p]$ with characteristic different from $p$)
 \cite{Ro,GM4}, the modular representations
 $F(\Si_g(2c))$ are no longer irreducible 
except if $g\leq 1$ or $(g,c)=(2,0)$.  
In all other cases, $F(\Si_g(2c))$ has a 
non-trivial  
 composition series with 
two irreducible factors
 \cite{GM4}. 
One can show \cite{M} that the representations of the
 mapping class group on these two factors induce well-defined 
 representations of the finite
 symplectic group $\Sp(2g,\BF_p)$. 
Following \cite{GM4}, let us denote the
 dimensions of these 
 representations  by $\fe_g^{(2c)}(p)$ and
 $\fo_g^{(2c)}(p)$. 
(The notation comes from the fact that the numbers $\fe_g^{(2c)}(p)$ and
 $\fo_g^{(2c)}(p)$ count a certain type of ``even'' and ``odd''
colorings of a graph.) 
Theorem~\ref{th1} has the
following corollary, which motivated this note. 

\begin{cor}\lbl{cor1} For $p=5$ and $g\geq 4$, we have the following equality of
  dimensions 
of irreducible 
 $\Sp(2g,\BF_p)$-representations: 
$$
      (\dim V(g,g  -m))_{m=0,1,2,3} 
= (  \fe_g^{(0)}(5),\fe_g^{(2)}(5),\fo_g^{(2)}(5),\fo_g^{(0)}(5))~.$$
This also holds  for $g<4$ if we
 interpret the left hand side as given by the classical formula
 (\ref{funk}) with $k=g-m$ 
in that case.
\end{cor}

At present, we don't  
know whether these equalities of dimensions
come from isomorphisms of the corresponding  
representations. We remark that for $g=1$ we have
identified the modular representation $F(\Sigma_1(2c))$ for any $p$ and
$c$ in \cite{GM2}:
in this case, one has $\fo_1^{(2c)}(p)=0$ and $\fe_1^{(2c)}(p)=
(p-1)/2-c$, and $F(\Sigma_1(2c))$ is isomorphic to the space of homogeneous
polynomials over $\BF_p$ in two variables of total degree $(p-3)/2-c$. It is
well-known that 
these representations
of $\Sp(2,\BF_p)=\SL(2,\BF_p)$ 
are 
the unique
(up to isomorphism) 
irreducible  
representations with these dimensions.

\vskip 8pt

The proof of Corollary~\ref{cor1} is based on explicit formulas for
the dimensions  $\fe_g^{(2c)}(p)$ and $ \fo_g^{(2c)}(p)$ which are
similar to the expression in
Theorem~\ref{th1}. These formulas are obtained as follows. First of
all, the 
 sum
\begin{equation} \lbl{sum1} D_g^{(2c)}(p)= \fe_g^{(2c)}(p)+\fo_g^{(2c)}(p)
\end{equation} is the dimension of  $F(\Si_g(2c))$, which is the same
as the rank of $\BS(\Si_g(2c))$, and thus the same as the dimension of
the $\SO(3)$-TQFT vector space 
$\BS(\Si_g(2c))\otimes {\mathbb C}$.
This
dimension is
given by the celebrated Verlinde formula
\begin{equation}\lbl{Verl} 
D_g^{(2c)}(p)=\left( \frac p 4\right)^{g-1} \sum_{j=1}^{(p-1)/2}
\left(\sin \frac { \pi j(2c+1)}{p}\right) \left( \sin \frac { \pi j}{p} \right)^{1-2g}
.\end{equation}   
We  
like to think of our expression for the dimensions of Gow's irreducible representations in
Theorem~\ref{th1} as an analogue of the Verlinde formula.

Next, we showed in \cite{GM4} that 
the difference 
\begin{equation} \lbl{sum2} \delta_g^{(2c)}(p)=
  \fe_g^{(2c)}(p)-\fo_g^{(2c)}(p) 
\end{equation} 
is given by  
the following  
formula:
\begin{equation}\lbl{d-Verl}  \delta_g^{(2c)}(p)= (-1)^c \,\frac{4}{p}
  \sum_{j=1}^{(p-1)/2}  \left(\sin\frac{\pi j (2c+1)}{p}\right) \left(\sin\frac{\pi j}{p}\right) \Lambda_j^g
\end{equation} 
where  $$\Lambda_j= \lceil (p-1) /4 \rceil + 2 \sum_{k=1}^{(p-3)/2}(-1)^k \lceil
(p-2k-1)/ 4 \rceil   \cos (2kj\pi/p)$$ 

Here, $\lceil x \rceil$ is the smallest integer $\geq x$.  
It is 
now 
easy to convert formulas 
(\ref{Verl}) and (\ref{d-Verl}) into
 explicit 
formulas for the 
dimensions $\fe_g^{(2c)}(p)$ and $ \fo_g^{(2c)}(p)$ of our irreducible  
$\Sp(2g,\BF_p)$-representations 
coming from  TQFT.

\vskip 8pt \noindent{\bf Remark.} We used these formulas in \cite{GM4} to show that the dimensions $\fe_g^{(2c)}(p)$ and $ \fo_g^{(2c)}(p)$ are
given by (inhomogeneous) polynomials in $p$ and $c$ of total degree
 $3g-2$, whose leading coefficients can be expressed in terms of Bernoulli
 numbers. 

\begin{proof}[Proof of Theorem~\ref{th1}] Let $q$ be a primitive
  $p$-th root of unity. Let $S=(S_{ij})$ be the
$(p-1)\times (p-1)$  matrix
  defined by
 \begin{equation*} \qquad \qquad \qquad \qquad S_{ij}=(-1)^{ij} (q^{ij}-q^{-ij}) \qquad \quad (i,j=1,2, \ldots, p-1)~.\end{equation*}
Let
   $Q=\diag(2+(-1)^j(q^j+q^{-j}))_{j=1,2,,\ldots,p-1}$, and let  $\I$
   denote the identity
  matrix. 

\begin{lem}\lbl{lem4} We have 
\begin{align}\lbl{MZ1} S^2&=- {2p} \,\I\\
\lbl{MZ} 
C&= S Q S^{-1}= - \textstyle {\frac{1}{2p}} S Q S 
\end{align} 
\end{lem}
\vskip 8pt

\noindent{\bf Remark.} 
This lemma (or some variants of it) is well-known in TQFT,
where it is used to
    prove the Verlinde formula (\ref{Verl}). In fact, our matrix $S$ is
    (up to some rescaling) the $S$-matrix of  the $\SU(2)$-TQFT at
    level  $p-2$. This TQFT corresponds to the $V_{2p}$-theory of
    \cite{BHMV2}. In the notation of that paper, the matrix $C$ is the
    matrix of multiplication  by $2+z$ in the standard basis of the
    Verlinde algebra $V_{2p}(\mathrm{torus})$, and the lemma
    follows from the diagonalization procedure on page 913 of
    \cite{BHMV2}. However, it is 
also
easy to prove the lemma by a direct
    computation, which  we include for the convenience of the
    reader.  

\begin{proof}[Proof of Lemma~\ref{lem4}] We have 
\begin{align} \notag (S^2)_{ij}&=\sum_{k=1}^{p-1} S_{ik}S_{kj}=
  \sum_{k=1}^{p-1} (-1)^{(i+j)k} (q^{ik}-q^{-ik}) (q^{kj}-q^{-kj})\\
\lbl{eq3} &=\sum_{k=1}^{p-1} (-1)^{(i+j)k}(q^{(i+j)k} + q^{-(i+j)k} -
q^{(i-j)k}- q^{-(i-j)k})
\end{align} This implies that the diagonal terms of $S^2$ are $(S^2)_{ii}=-2p$,
since 
\begin{equation} \lbl{eq2}\sum_{k=1}^{p-1} q^{2ik}=-1 \hspace{2cm} (1\leq
  i\leq p-1).
\end{equation} It remains to show that $S_{ij}=0$ if $i\neq j$. If
$i+j$ is even, this follows again from (\ref{eq2}), while if $i+j$ is
odd, this follows from

\begin{equation} \lbl{eq1} \sum_{k=1}^{p-1}(-1)^k (q^{nk} +
q^{-nk})=0 \hspace{2cm} (n\in \BZ)\end{equation} (To prove (\ref{eq1}), observe that the
$k$-th summand cancels the $(p-k)$-th summand.) Thus we have proved
the first statement of the lemma.  To prove the second statement
(\ref{MZ}), it suffices to 
observe that 
\begin{align*} (CS)_{ij}&=\sum_{k=1}^{p-1} C_{ik}S_{kj}= 2 S_{ij} +
  S_{i+1,j} + S_{i-1,j}\\
&= \left(2+(-1)^j(q^j+q^{-j})\right)S_{ij}\\
&=(SQ)_{ij}
\end{align*} (This computation is also valid for $i=1$ and $i=p-1$ since
$S_{0j}=0=S_{pj}$.) Thus we have $CS=SQ$, which implies (\ref{MZ}).
\end{proof}

Using the lemma, we now have 
 \begin{align*}(C^g&)_{n,1}=-
  \frac{1}{2p} \sum_{j=1}^{p-1} S_{nj}\Big(2+(-1)^j(q^{j}+q^{-j})\Big)^g S_{j1}\\
&=  -
  \frac{1}{2p}
  \sum_{j=1}^{p-1}(-1)^{(n+1)j}(q^{nj}-q^{-nj})(q^{j}-q^{-j})\Big(2+(-1)^j(q^{j}+q^{-j})\Big)^g\\
=&  -
  \frac{1}{2p}
  \sum_{j=1}^{(p-1)/2}(q^{nj}-q^{-nj})(q^{j}-q^{-j}) 
\Big( (2+q^{j}+q^{-j})^g  - (-1)^{n} (2-q^{j}-q^{-j})^g   \Big) 
\end{align*}
(In the last step, we have grouped together the term with $j$ and the
term with $p-j$. It helps to consider separately the cases  where $j$
is even and where $j$ is odd.)  Substituting now $q=e^{2\pi i/p}$ into this formula gives
Theorem~\ref{th1}.
\end{proof}
\begin{proof}[Proof of Corollary~\ref{cor1}]  Let $G_j$ be the Galois
  automorphism of the cyclotomic field $\BQ(q)$ defined by
  $G_j(q)=q^j$. Formulas (\ref{Verl}) and (\ref{d-Verl}) can be
  expressed as follows \cite{GM4}: 
\begin{align}\lbl{V1}D_g^{(2c)}&=-
  \frac{1}{p} \sum_{j=1}^{(p-1)/2}  G_{j}
  \Bigg((q^{2c+1}-q^{-2c-1})(q-q^{-1})\Big(\frac{-p}{(q-q^{-1})^2}\Big)^g\Bigg)\\ \lbl{V2}
\delta_g^{(2c)}&=-
  \frac{1}{p} \sum_{j=1}^{(p-1)/2}  G_{j}
  \Bigg((q^{2c+1}-q^{-2c-1})(q-q^{-1})(-1)^c \Lambda^g\Bigg)
\end{align} where $$\Lambda = \lceil (p-1) /4 \rceil +  \sum_{k=1}^{(p-3)/2}(-1)^k \lceil
(p-2k-1)/ 4 \rceil  (q^{2k}+q^{-2k})~.$$ 

Now let $p=5$ so that $1+q+q^2+q^3+q^4=0$. Thus
$$\Lambda=1-q^2-q^{-2}=2+q+q^{-1}$$ and one checks that
$(q-q^{-1})^2(2-q-q^{-1})=-5$ so that $$ \frac{-5}{(q-q^{-1})^2}=
2-q-q^{-1}~.$$
Therefore \begin{align*}(C^g)_{n,1}&=\frac{-1}{10}\sum_{j=1}^{2} G_j\left(
  (q^{n}-q^{-n})(q-q^{-1 })\Bigg( \Lambda^{g} -(-1)^{n} \left(\frac{-5}{(q-q^{-1})^2}\right)^{g}
  \Bigg)\right)\\
\end{align*}
Comparing this with (\ref{V1}) and (\ref{V2}) proves the corollary.
\end{proof}
\vskip 8pt

\noindent{\bf Remark.}  We see that for $p\geq 7$, there is no obvious
    relation between  $2+q+q^{-1}$ and
    $2-q-q^{-1}$ on the one hand, and $\Lambda$ and
    $-p/((q-q^{-1})^2)$ on the other hand. This explains why for
    $p\geq 7$ we found
    no  relation between the dimensions of 
Gow's representations  
$V(g,k)$ on the one hand,
    and the dimensions $\fe_g^{(2c)}(p)$ and
 $\fo_g^{(2c)}(p)$  
of our irreducible $\Sp(2g,\BF_p)$-representations coming from TQFT  
on the other hand.  
It is, however, intriguing that we get as many irreducible
representations as Gow (namely $p-1$, since we have two of them for
each $0\leq c\leq 
(p-3)/2$). We wonder about how to characterize our
representations among the irreducible representations of
$\Sp(2g,\BF_p)$ in characteristic $p$.

\vskip 8pt
\noindent{\em Acknowledgments.} We thank Marc Cabanes, Jens Carsten
Jantzen, and especially Henning Haahr Andersen, for helpful discussions about highest weight
representations in positive characteristic.

\end{document}